\title{ ~~\\ The formal series Witt transform}
\author{Pieter Moree}
\documentclass[12pt]{article}
\usepackage{amssymb, latexsym, amsfonts}
\def\@ptsize{2}
\newtheorem{Thm}{Theorem}

\newtheorem{Lem}{Lemma}
\newtheorem{Def}{Definition}

\newtheorem{Prop}{Proposition}
\newcommand{\qed}{\hfill $\Box$}

\begin{document}
\date{}
\maketitle
{\def\thefootnote{}
\footnote{{\it Mathematics Subject Classification (2000)}.
Primary 05A19; Secondary 11B75, 17B01}}
\begin{abstract}
\noindent Given a formal power series $f(z)\in \mathbb C[[z]]$ we
define, for any positive integer $r$, its $r$th Witt
transform, ${\cal W}^{(r)}_f$, by
${\cal W}^{(r)}_f(z)={1\over r}\sum_{d|r}\mu(d)f(z^d)^{r/d}$, where $\mu$
denotes
the M\"obius function. The Witt transform generalizes the necklace
polynomials, $M(\alpha;n)$, that occur in the cyclotomic identity
$${1\over 1-\alpha y}=\prod_{n=1}^{\infty}(1-y^n)^{-M(\alpha;n)}.$$
Several properties of ${\cal W}^{(r)}_f$ are established.
Some examples relevant to number theory are considered.
\end{abstract}
\section{Introduction}
The polynomial
$$M(\alpha;n)={1\over n}\sum_{d|n}\mu({n\over d})\alpha^d$$
is called the {\it necklace polynomial} in \cite{MR} and arises naturally in
many combinatorial problems.
This polynomial is of degree $n$ in $\alpha$ with rational coefficients
and takes on integer
values for integer arguments 
(it is a so-called {\it integral polynomial}). 
Taking $n$ to be a prime we infer from this 
Fermat's little theorem and indeed, in this context the 
polynomials $M(\alpha;n)$ were first studied (starting with Gauss), 
see \cite[Chapter 11]{D}. If $\alpha\in \mathbb Z_{\ge 1}$, $n\ge 3$ 
and $n\not\equiv 2({\rm mod~}4)$, then $M(\alpha;n)$ is even \cite{CM}. 
It is not difficult to show
\cite[Lemma 3]{Moreeconstant} that for $n\ge 1$ and $\alpha>1$, with
$\alpha$ real,
$M(\alpha;n)>0$. 
A sequence $\{a_n\}_{n=1}^{\infty}$ of non-negative integers is said to be 
{\it exactly realizable} if there is a set $X$ and a map $T:X\rightarrow X$ for
which $\#\{x\in X|T^nx=x\}=a_n$ for all $n\ge 1$. 
Puri and Ward \cite{PW} proved that a sequence $\{a_n\}_{n=1}^{\infty}$ of
non-negative integers is exactly realizable iff $\sum_{d|n}\mu(n/d)a_d$ is
non-negative and divisible by $n$ for all $n\ge 1$.
Taking $X=\mathbb C$ and
$T$ the map that sends $z$ to $z^{\alpha}$, we see that the sequence 
$\{{\alpha}^n\}_{n=1}^{\infty}$ 
is exactly realizable for $\alpha\in \mathbb Z_{\ge 1}$ and hence, by
the result of Puri and Ward, that  $\{M(\alpha;n)\}_{n=1}^{\infty}$ consists 
of non-negative integers only.\\
\indent The necklace polynomial $M(\alpha;n)$ got its name since
it can be
interpreted as enumerating non-periodic circular strings of $n$ beads that
can be strung from beads
of at most $\alpha$ distinct colours. It is called the Witt formula when
used to count the number
of monic irreducible polynomials of degree $n$  over the finite field
$\mathbb F_q$, with $q$ a prime power.
It also gives the dimension of the subspace spanned by the homogeneous
elements of degree $n$ in the
free Lie algebra over a set of $\alpha$ elements (this is
the original context in which Witt \cite{W} discovered his formula).
The necklace polynomial also arises in the context of
Philip Hall's commutator collecting algorithm, 
see e.g. \cite[Chapter 11]{Hall}. Golomb \cite{G} showed that the maximum number of words
possible in a bounded synchronization delay code with word length $n$ over an alphabet
of $\alpha$ elements equals $M(\alpha;n)$.   
More recently necklace
polynomials also arose in the study of multiple zeta series \cite{Ho}. There are also connections 
with the theory of formal groups \cite{L}.\\
\indent The
cyclotomic identity states, that as formal series
we have
$${1\over 1-\alpha y}=
\prod_{n=1}^{\infty}\left({1\over 1-y^n}\right)^{M(\alpha;n)}.$$
Using logarithmic differentation and M\"obius inversion the
cyclotomic identity is easily established. Metropolis and Rota
\cite{MR}
gave the first natural, i.e. bijective proof of the cyclotomic
identity, that is a proof which is entirely set-theoretic,
where set-theoretic constructions are made to correspond
biuniquely to the algebraic operations on formal power series.
They also noted the following properties of
the necklace polynomials (where $(a,b)$ denotes the greatest
common divisor and $[a,b]$ the least common multiple of
the integers $a$ and $b$).
\begin{Thm}
\label{metro}
{\rm (Metropolis and Rota {\rm \cite{MR}})}.\\
{\rm 1)} We have, for any positive integers $\alpha,\beta$ and $n$:
$$M(\alpha \beta; n)=\sum_{[i,j]=n}(i,j)M(\alpha;i)M(\beta;j),$$
where the sum ranges over all positive integers $i$ and $j$ with
$[i,j]=n$.\\
{\rm 2)} We have, for any positive integers $\beta,r$ and $n$:
$$M(\beta^r;n)=\sum_{[j,r]=nr}{j\over n}M(\beta;j),$$
where the sum ranges over the integers $j$ with $[r,j]=nr$.
\end{Thm}
Metropolis and Rota write regarding the above identities: `We shall
be concerned with some remarkable identities satisfied by the polynomials
$M(\alpha,n)$, which apparently have not been previously noticed'.
However, certainly part 1 of their result was known long 
before, see e.g. \cite{C}, but Metropolis and Rota
gave the first combinatorial proof.\\
\indent In this paper we consider a generalization of the cyclotomic identity and
hence
also of the necklace polynomial:
\begin{Def}
For $f(z)\in \mathbb C[[z]]$ and $r\ge 1$ any integer, let
$${\cal W}^{(r)}_f(z)={1\over
r}\sum_{d|r}\mu(d)f(z^d)^{r/d}=\sum_{j=0}^{\infty}m_f(j,r)z^j.$$
\end{Def}
With this definition the cyclotomic
identity generalizes as follows:
\begin{Thm} {\rm (Moree {\rm \cite{PNew}})}.
\label{een}
Suppose that $f(z)\in \mathbb Z[[z]]$.
Then, as formal power series in $y$ and $z$, we have
$${1\over
1-yf(z)}=\prod_{j=0}^{\infty}\prod_{k=1}^{\infty}(1-z^jy^k)^{-m_f(j,k)}.$$
\end{Thm}
Note that if we take $f(z)=\alpha$, the cyclotomic identity is obtained
and that
${\cal W}^{(r)}_f(z)=M(\alpha;r)$.\\
\indent The following general properties of ${\cal W}^{(r)}_f$ will be
established in this note. Parts 4 and 5 generalize Theorem \ref{metro}.
\begin{Thm}
\label{main1}
Let $r\ge 1$ be an integer. Let $f,g\in \mathbb C[[z]]$.
\item{\rm 1)} We have ${\cal W}_{z^kf}^{(r)}(z)=z^{kr}{\cal W}^{(r)}_f(z)$.
\item{\rm 2)} We have $\sum_{d|r}{r\over d}{\cal W}^{(r/d)}_f(z^d)=f(z)^r$,
\item{\rm 3)} We have
$$(-1)^r {\cal W}^{(r)}_{-f}(z)=\cases{{\cal W}^{(r)}_f(z)+{\cal
W}^{(r/2)}_f(z^2) &if $r\equiv 2({\rm mod~}4)$;\cr
{\cal W}^{(r)}_f(z) &otherwise.}$$
\item{\rm 4)} We have 
$${{\cal W}_{fg}}^{(r)}(z)=
\sum_{[i,j]=r}(i,j){\cal W}_f^{(i)}(z^{r\over i}){\cal W}_g^{(j)}(z^{r\over j}),$$
where the sum is over all positive integers $i$ and $j$ with $[i,j]=r$.
\item{\rm 5)} We have
$${\cal W}_{f^k}^{(r)}(z)=\sum_{[j,k]=rk}{j\over r}{\cal W}_f^{(j)}(z^{rk\over j}),$$
where the sum is over all positive integers $j$ with $[j,k]=rk$.
\item{\rm 6)} Let $v$ and $w$ be positive integers. Then
$${\cal W}_{f^{w/(v,w)}g^{v/(v,w)}}^{(r)}(z)=\sum 
{(vi,wj)\over (v,w)}{\cal W}_f^{(i)}(z^{r\over i})
{\cal W}_g^{(j)}(z^{r\over j}),$$
where the sum ranges over the set $\{i,j:~{ij\over (vi,wj)}={r\over (v,w)}\}$.
\end{Thm}
Tbe latter three properties simplify if one puts ${\cal C}_f^{(r)}(z)=
r{\cal W}_f^{(r)}(z)$; the above identities then hold 
with ${\cal W}$ replaced by ${\cal C}$ and $(i,j)$ (in part 4), $j/r$ 
(in part 5) and $(vi,wj)/(v,w)$ (in part 6) left out.\\
\indent In the following result the coefficients of $f$ are assumed to be integers
(recall that a polynomial $f$ is self-reciprocal if $z^{{\rm
deg}f}f(1/z)=f(z)$).
\begin{Thm}
\label{main2}
\item{\rm 1)} If  $f(z)\in \mathbb Z[z]$ is self-reciprocal, then so is ${\cal
W}_f^{(r)}(z)$.
\item{\rm 2)} If $f(z)\in \mathbb Z[[z]]$, then ${\cal W}^{(r)}_f(z)\in \mathbb
Z[[z]]$.
\item{\rm 3)} If $f(z)\in \mathbb Z_{\ge 0}[[z]]$, then ${\cal W}^{(r)}_f(z)\in
\mathbb Z_{\ge 0}[[z]]$.
\item{\rm 4)} If $f(z)\in \mathbb Z_{\le 0}[[z]]$, then $(-1)^r{\cal
W}^{(r)}_f(z)\in \mathbb Z_{\ge 0}[[z]]$.
\item{\rm 5)} If $f(z)\in \mathbb Z_{\ge 0}[[z]]$ and $g(z)-f(z)\in \mathbb
Z_{\ge 0}[[z]]$, then
${\cal W}^{(r)}_g(z)-{\cal W}^{(r)}_f(z)\in \mathbb Z_{\ge 0}[[z]]$.
\end{Thm}
The final and (deepest) result is concerned with the monotonicity of the coefficients
of ${\cal W}_f^{(r)}(z)$.
\begin{Thm}
\label{main3}
Let $f(z)=\sum_j a_jz^j\in \mathbb Z_{\ge 0}[[z]]$. In parts {\rm 1} and 
{\rm 2} it is assumed
that $a_0>0$. In the remaining parts it is assumed in addition that 
$\{a_j\}_{j=0}^{\infty}$ is a non-decreasing sequence.
\item{\rm 1)} Let $k\ge 2$. The
sequence $\{m_f(k,r)\}_{r=1}^{\infty}$ is non-decreasing.
\item{\rm 2)} Let $k\ge 3$. The
sequence $\{(-1)^rm_{-f}(k,r)\}_{r=1}^{\infty}$ is non-decreasing.
\item{\rm 3a)} If $r\ge 1$, then $m_f(k,r)\ge 1$.
\item{\rm 3b)} If $r\ge 1$, then the
sequence $\{m_f(k,r)\}_{k=0}^{\infty}$ is non-decreasing. 
\item{\rm 3c)} If $r\ge 3$, the 
sequence $\{m_f(k,r)\}_{k=2}^{\infty}$ is strictly increasing.
\item{\rm 4a)} If $r\ge 1$, then $(-1)^rm_{-f}(k,r)\ge 1$.
\item{\rm 4b)} If $r\ge 1$, then the
sequence $\{(-1)^rm_{-f}(k,r)\}_{k=0}^{\infty}$ is non-decreasing. 
\item{\rm 4c)} If $r\ge 3$, the 
sequence $\{(-1)^rm_{-f}(k,r)\}_{k=2}^{\infty}$ is strictly increasing.
\end{Thm}
The condition that $\{a_j\}_{j=0}^{\infty}$ be non-decreasing for parts 3 and 4
seems to be rather stringent, but actually cannot be dropped.\\
\indent Note that $m_f(0,r)=M(f(0);r)$. The monotonicity 
aspects of necklace polynomials are not covered by Theorem \ref{main3}, but
are easily determined using the same methods:
\begin{Prop}
\label{wassendwater}
Let $\beta(2)=\beta(3)=2$ and $\beta(k)=1$ for $k\ge 4$. Let $c\ge 2$ be an integer.
The sequence $\{M(c;r)\}_{r=\beta(c)}^{\infty}$ is strictly increasing. Let $r\ge 1$.
The sequence $\{M(c;r)\}_{c=1}^{\infty}$ is strictly increasing.
\end{Prop}

\indent The remaining part of the paper is concerned with the proof of Theorems
\ref{main1}, \ref{main2}, \ref{main3} and the latter proposition.
In the next section some lemmas on circular words are established.
In the final section some examples are discussed.

\section{Circular words and Witt's dimension formula}
We will make use of an easy result on cyclic words. A word $a_1\cdots a_n$ is
called {\it circular} 
or {\it cyclic} if $a_1$ is regarded as following $a_n$, where
$a_1a_2\cdots a_n$,
$a_2\cdots a_na_1$
and all other cyclic shifts (rotations)
of $a_1a_2\cdots a_n$ are regarded as the same word. A circular
word of length $n$ may conceivably be given by repeating a segment of
$d$ letters $n/d$
times, with $d$ a divisor of $n$. Then we say the word is of {\it period} $d$.
Each word belongs
to an unique smallest period; the {\it minimal period}.\\
\indent Consider
circular words of length $n$ on an alphabet $x_1,\dots,x_r$ consisting of
$r$ letters.
The total number of ordinary words such that $x_i$ occurs $n_i$ times equals
$({n\atop n_1,\dots,n_r})$, where $n_1+\dots+n_r=n$ and
$$\left({n\atop n_1,\dots,n_r}\right)={n!\over n_1!\cdots n_r!}.$$
Let $M(n_1,\dots,n_r)$ denote
the number of circular words of length $n_1+\dots+n_r=n$ and
minimal period $n$ (often called {\it aperiodic words}) such that the letter $x_i$ appears
exactly $n_i$ times.
This leads to the formula
\begin{equation}
\label{witty}
\left({n\atop n_1,\dots,n_r}\right)=\sum_{d|{\rm
gcd}(n_1,\dots,n_r)}{n\over d}M({n_1\over d},{n_2\over d},
\dots,{n_r\over d}).
\end{equation}
whence it follows by M\"obius inversion that
\begin{equation}
\label{basaal}
M(n_1,\dots,n_r)={1\over n}
\sum_{d|{\rm gcd}(n_1,\dots,n_r)}\mu(d)\left({{n\over d}\atop {n_1\over
d},\cdots,{n_r\over d}}\right).
\end{equation}
Note that $M(n_1,\dots,n_r)$ is totally symmetric in the variables $n_1,\dots,n_r$.
If $\{n_j\}_{j=1}^{\infty}$ is a sequence such that
there is a $k$ for which $n_j=0$ for every $j\ge k+1$, then we define
$M(\{n_j\}_{j=1}^{\infty})=M(n_1,\dots,n_k)$.
Note that
\begin{equation}
\label{handig}
{1\over n}\sum_{d|n}\mu(d)(z_1^d+\cdots+z_r^d)^{n\over
d}=\sum_{n_1+\cdots+n_r=n\atop n_j\ge 0}M(n_1,\cdots,n_r)z_1^{n_1}
\cdots z_r^{n_r}.
\end{equation}
\indent It turns out that the numbers $M(n_1,\dots,n_r)$ are related to  counting
so called
basic commutators in group theory \cite[Chapter XI]{Hall}. These numbers
also occur in a classical result
in Lie theory, namely Witt's formula for the homogeneous subspaces of a
finitely generated
free Lie algebra $L$: if $H$ is the subspace of $L$ generated by all
homogeneous elements of
multidegree $(n_1,\dots,n_r)$, then dim$(H)=M(n_1,\dots,n_r)$, where
$n=n_1+\dots+n_r$. In the Lie algebra context the cyclotomic identity is 
interpreted as a denominator identity related to the free Lie algebra, see
e.g. \cite{KK2, KK}. As the referee pointed out the symmetric polynomials in
(\ref{handig}) have been studied in the theory of symmetric functions and
have applications to counting permutations with certain properties, for 
example unimodal permutations, see e.g. Thibon \cite{Thi}, and
permutations with prescribed descent set, see e.g. Gessel
and Reutenauer \cite{GR}.\\
\indent Using (\ref{witty}) one infers (on taking the logarithm of either side and
expanding it as a formal series) that
$${1\over 1-z_1-\cdots
-z_r}=\prod_{n_1,\cdots,n_r=0}^{\infty}(1-z_1^{n_1}\cdots
z_r^{n_r})^{-M(n_1,\cdots,n_r)},$$
where $(n_1,\cdots,n_r)=(0,\cdots,0)$ is excluded in the product. 
It is a consequence of the latter identity, that if $a_j\in \mathbb
Z_{\ge 0}$ for $1\le j\le n$, then
as formal series we have $1-a_1t-\cdots -a_nt^n=\prod_{k\ge
1}(1-t^k)^{e_k}$, with $e_k\in \mathbb Z_{\ge 0}$. This
was proved earlier in \cite[Lemma 4]{Moreeconstant} using zeta functions
of finite automata.\\

\subsection{A sign twisted variation of $M(n_1,\dots,n_r)$}
In our considerations a sign twisted variation of (\ref{basaal}) comes up.
\begin{Lem}
\label{minnetje}
Let $k$ and $r$ be positive integers with $k\le r$ and $n_j$ non-negative
integers and put
$t_k=n_1+\dots+n_k$ and $n=n_1+\dots+n_r$. Put
$$V_k(n_1,\dots,n_r)={(-1)^{t_k}\over
n}\sum_{d|{\rm gcd}(n_1,\dots,n_r)}\mu(d)(-1)^{t_k\over d}
\left({{n\over d}\atop {n_1\over d},\cdots,{n_r\over d}}\right)$$
Then $$V_k(n_1,\dots,n_r)=$$
$$\cases{
M(n_1,\dots,n_r)+M({n_1\over 2},\dots,{n_r\over 2}) &if $t_k\equiv 2({\rm
mod~}4)$ and
$2|{\rm gcd}(n_1,\dots,n_r)$;\cr
M(n_1,\dots,n_r) & otherwise.}$$
\end{Lem}
{\it Proof}. The only not immediately obvious case is when  $t_k\equiv
2({\rm mod~}4)$ and
$2|{\rm gcd}(n_1,\dots,n_r)$. So assume we are in this case. Note that then
at least one of the $n_j$ is congruent to $2({\rm mod~}4)$. Write
$M(n_1,\dots,n_r)$ as
$S_{\rm odd}+S_{\rm even}$, where in $S_{\rm odd}$ all terms
with $d$ odd are collected.
Thus $M(n_1,\dots,n_r)=S_{\rm odd}+S_{\rm even}$. We have
$$V_k(n_1,\dots,n_r)=S_{\rm odd}-S_{\rm
even}=M(n_1,\dots,n_r)-2S_{\rm even}.$$
Using that at least one of the
$n_j$ satisfies $n_j\equiv 2({\rm mod~}4)$, we infer that $2S_{\rm
even}=-M(n_1/2,\dots,n_r/2)$. \qed\\

\noindent Note that
$$
{1\over
n}\sum_{d|n}\mu(d)(-z_1^d-\cdots-z_k^d+z_{k+1}^d+\cdots+z_r^d)^{n\over d}=$$
\begin{equation}
\label{handig2}
\sum_{n_1+\cdots+n_r=n\atop n_j\ge 0}(-1)^{n_1+\dots+n_k}V_k(n_1,\dots,n_r)z_1^{n_1}
\cdots z_r^{n_r}.
\end{equation}

\noindent {\tt Remark 1}. The numbers $M(k,m-k)$ and $V_1(k,m-k)$ were
already studied in a different guise
around 1900 by R. Daublebsky von Sterneck
(see e.g. \cite[Vol. II, pp. 222-264]{B}).
Daublebsky von Sterneck showed that the number of ways of selecting $k$ parts,
respectively $k$ distinct parts, from $0,1,\dots,m-1$
so that their sum is congruent to $1({\rm mod~}m)$ equals $M(k,m-k)$,
respectively $V_1(k,m-k)$. Simple proofs of Daublebsky von Sterneck's results
were later given by Ramanathan \cite{R2}, see also \cite{BFW, EJP}. Ramanathan 
uses properties of Ramanujan sums and in \cite{BFW,EJP} the authors make
use of Gauss polynomials. 
Let $\varphi$ denote Euler's totient function.
The function $\Phi(k,n)=\varphi(n)\mu(n/(k,n))/\varphi(n/(k,n)$
(called von Sterneck function by some
authors) was introduced by Daublebsky von Sterneck in this context. He proved several
of its properties. The von Sterneck function, however, is equal to the Ramanujan
sum $c_n(k)$, which was introduced later by Ramanujan.\\

\noindent {\tt Remark 2}. The numbers $M(n_1,\cdots,n_r)$ were interpreted as dimensions
by Witt (see the previous section). The numbers $V_k(n_1,\dots,n_r)$
can also be interpreted as dimensions (in the context of
free Lie superalgebras), see \cite{P}.\\

\noindent {\tt Remark 3}. By setting $V_0(n_1,\dots,n_r,m)=M(n_1,\dots,n_r)$ it is
possible to deal with $M(n_1,\dots,n_r)$ and $V_k(n_1,\dots,n_r)$ in a more uniform
way. For reasons of exposition this route has not been chosen.

\subsection{Lyndon words}
If $w$ is a circular word counted by $M(n_1,\dots,n_r,m)$ we can choose amongst
the rotations of $w$ one which is lowest with respect to a given lexicographical
order (since we work with numbers as letters it is most natural to say
that $i<j$ if $i<j$ as natural numbers). This is the idea of Lyndon words, which
we now describe more precisely.\\
\indent If $\cal A$ is an alphabet (assumed to be finite for simplicity of description), 
let ${\cal A}^*$ be the set
of words with letters from ${\cal A}$ and ${\cal A}^+$ the set of non-empty words.
Suppose we have a total order on ${\cal A}$. We extend the total ordering to ${\cal A}^+$ in
the following way: For any $u,v\in {\cal A}^+$, $u<v$ iff either $v\in u{\cal A}^+$
or $u=ras$, $v=rbt$, with $a<b$; $a,b\in {\cal A}$; $r,s,t\in {\cal A}^*$. By
definition a 
{\it Lyndon word} is an aperiodic word that is minimal amongst all the rotations
of it. E.g. for ${\cal A}=\{a,b\}$ and $a<b$, the list of first Lyndon words is
$\{a,b,ab,aab,abb,aaab,aabb,\cdots\}$. Let $L$ denote the set of Lyndon words.
The following proposition is quite useful.
\begin{Prop}
\label{Lyndon}
A word $w\in {\cal A}^+$ is a Lyndon word iff $w\in {\cal A}$ or $w=lm$ with
$l,m\in L$, $l<m$.   
\end{Prop}
{\it Proof}. Cf. the proof of Proposition 5.1.3 of \cite{Lot}. \qed\\

\noindent The above proposition shows that given a Lyndon word $w$, the word $wz$
is also Lyndon (unless $w=z$), where $z$ is the letter which is
highest in the total order on ${\cal A}$. We call this procedure {\it Lyndon extension}.\\
\indent In \cite[Section 5.3]{Lot} the connection of Lyndon words
with free Lie algebras is discussed, cf. \cite{V}.

\subsection{Monotonicity}
In this section monotonicity properties of $M(n_1,\cdots,n_r)$ and $V_k(n_1,\dots,n_r)$ 
are being considered.
 
\begin{Prop} $~$\\
{\rm 1)} Suppose that $n_2\ge 1$, then $M(0,n_2,\dots,n_r)\le
M(1,n_2-1,n_3,\dots,n_r).$\\
{\rm 2)} If $n_2\ge 2$, then $M(0,n_2,\dots,n_r)<M(1,n_2-1,n_3,\dots,n_r)$.\\
Let $1\le k\le r$.\\
{\rm 3)} Suppose that $n_2\ge 1$, then $V_k(0,n_2,\dots,n_r)\le
V_k(1,n_2-1,n_3,\dots,n_r).$\\ 
{\rm 4)} If $n_2\ge 2$ 
and $n_2+\cdots+n_r\ge 3$, then $V_k(0,n_2,\dots,n_r)<V_k(1,n_2-1,n_3,\dots,n_r)$.
\end{Prop}
{\it Proof}. 1) Let $w$ be a Lyndon word counted by $M(0,n_2,\dots,n_r)$. It
starts with a 2 and does not contain a 1. Replace this 2 by a 1. This
yields a Lyndon word counted by $M(1,n_2-1,n_3,\dots,n_r)$. Since this
procedure is injective it follows that $M(0,n_2,\dots,n_r)\le
M(1,n_2-1,n_3,\dots,n_r)$.\\
2) Since by assumption $n_2\ge 2$, the set $M(1,n_2-1,n_3,\dots,n_r)$ counts
at least one Lyndon word of the form $1z2$. Since $2z2$ is not a Lyndon word, the
claimed inequality follows.\\
3) This follows from part 1, together with Lemma \ref{minnetje} except
for the case where $t_k\equiv 2({\rm mod~}4)$ and
$2|{\rm gcd}(n_1,\dots,n_r)$, in which case we have to show that
$$
M(0,n_2,\dots,n_r)+M(0,{n_2\over 2},\dots,{n_r\over 2})
\le M(1,n_2-1,n_3,\dots,n_r).
$$
In case the second quantity in the latter inequality equals zero, we are
done by part 1, so assume that $M(0,n_1/2,\dots,n_r/2)\ge 1$. The Lyndon
words counted by $M(0,n_2,\dots,n_r)$ we deal with as before. If $2w$ is a Lyndon
word counted by $M(0,n_1/2,\dots,n_r/2)$, we consider the word
$1w2w$. Since $w$ does not contain a 1, it is a Lyndon word. It is counted
by $M(1,n_2-1,n_3,\dots,n_r)$ and $2w2w$ is not counted by 
$M(0,n_2,\dots,n_r)$.\\
4) This follows on using part 2 and noting in addition that $1w2w$ does not end in a 2.
 \qed\\

\noindent The idea of using Lyndon extension to
prove the next result was kindly communicated to
the author by Prof. F. Ruskey. Profs. Bryant \cite{Br} and
Petrogradsky \cite{Petro} proved part 1 of the next proposition using
Lie algebraic methods (Hall basis, respectively Lyndon-Shirshov basis of a free
Lie algebra).
\begin{Prop} \label{mooizeg} Let $1\le k\le r+1$.\\
\noindent {\rm 1)} Suppose that $n_1+\cdots+n_r\ge 1$.  Then
$\{M(n_1,\dots,n_r,m)\}_{m=0}^{\infty}$ is
a non-decreasing sequence. We have
$\{M(m,0,\dots,0)\}_{m=1}^{\infty}=\{1,0,0,\dots\}$ in the remaining case (i.e. 
the case $n_1+\cdots+n_r=0$).\\
{\rm 2)} Suppose that $n_1+\cdots+n_r\ge 1$.  Then
$\{V_k(n_1,\dots,n_r,m)\}_{m=0}^{\infty}$ is
a non-decreasing sequence. In the remaining case one has 
$\{V_k(m,0,\dots,0)\}_{m=1}^{\infty}=\{1,1,0,\dots\}$.
\end{Prop}
{\it Proof}. 1) First proof. This follows at once from Lyndon extension. Choose the 
Lyndon words as representatives of the circular words counted by
$M(n_1,\dots,n_r,m)$. Now concatenate each such word with
the letter $r+1$. By Proposition \ref{Lyndon} this yields another
Lyndon word which is counted by $M(n_1,\dots,n_r,m+1)$.
Since the concatenation is injective, it follows that
$M(n_1,\dots,n_r,m)\le M(n_1,\dots,n_r,m+1)$.\\
Second proof (by Dion Gijswijt). If
$M(n_1,\dots,n_r,m)=0$, there is nothing to prove,  so
assume that $M(n_1,\dots,n_r,m)\ge 1$.
If we have a circular word counted by $M(n_1,\dots,n_r,m)$ do
the following: if $m=0$ insert the letter $r+1$ anywhere in the sequence. This
yields an aperiodic circular word that is counted by $M(n_1,\dots,n_r,m+1)$. If
$m\ge 1$ look for a longest consecutive string of letters $r+1$ in
a circular word counted by $M(n_1,\dots,n_r,m)$ and insert
another letter $r+1$ after it. This clearly yields an aperiodic word counted by
$M(n_1,\dots,n_r,m+1)$. Since this extension procedure is injective, the result 
follows.\\
2) This follows from part 1, together with Lemma \ref{minnetje} except
for the case where $t_k\equiv 2({\rm mod~}4)$ and
$2|{\rm gcd}(n_1,\dots,n_r,m)$, in which case we have to show that
\begin{equation}
\label{geenidee}
M(n_1,\dots,n_r,m)+M({n_1\over 2},\dots,{n_r\over 2},{m\over 2})
\le M(n_1,\dots,n_r,m+1).
\end{equation}
In case the second quantity in the latter inequality equals zero, we are
done by part 1, so assume that $M(n_1/2,\dots,n_r/2,m/2)\ge 1$. The Lyndon
words counted by $M(n_1,\dots,n_r,m)$ we extend as before. 
Let $W$ be the set of words thus produced.
If $w$ is a Lyndon
word counted by $M(n_1/2,\dots,n_r/2,m/2)$, we consider the word
$wwz$, where $z$ stands for the letter $r+1$. By Proposition \ref{Lyndon} 
it follows that $wwz$ is a Lyndon word (note that $w\ne z$). This word is
counted by $M(n_1,\dots,n_r,w)$ and is not in $W$.\\
\indent For the final part of the assertion we use the easy observation that
$${(-1)^m\over m}\sum_{d|m}(-1)^{m\over d}\mu(d)=\cases{1 &if $m\le 2$;\cr
0 & otherwise.}$$
This concludes the proof. \qed\\

\noindent The following result sharpens Proposition \ref{mooizeg}.
\begin{Thm} 
\label{atlonglast}
Let $r\ge 1$ and $n_1,\dots,n_r$ be non-negative numbers.\\
{\rm 1)}  The
sequence $\{M(n_1,\dots,n_r,m)\}_{m=0}^{\infty}$ is strictly
increasing if $n_1+\cdots+n_r\ge 3$ or $r=2$ and $n_1=n_2=1$.\\
{\rm 2)} The sequence 
$\{V_k(n_1,\dots,n_r,m)\}_{m=0}^{\infty}$ is strictly
increasing if $n_1+\cdots+n_r\ge 3$ or $r=2$ and $n_1=n_2=1$, 
when $1\le k\le r+1$.
\end{Thm}
In the proof we make use of the following trivial result.
\begin{Lem}
\label{allebegin}
We have
$$M(0,m)=\cases{1 &if $m=1$;\cr 0 &if $m>1$,}{\rm ~and~}
M(1,m)=1{\rm~for~}m\ge 1.$$
Let $m\ge 0$. We have
$$M(2,m)=\cases{m/2 &if $m$ is even;\cr (m+1)/2 & otherwise.}$$ 
If gcd$(n_1,\dots,n_r)=1$ we have
$$M(n_1,\dots,n_r,m+1)={(n_1+\dots+n_r+m)\over m+1}M(n_1,\dots,n_r,m).$$
\end{Lem}

\noindent {\it Proof of Theorem} \ref{atlonglast}. 1)
If $n_1,n_2\ge 1$, then $1^{n_1}\cdots r^{n_r}$ is a Lyndon word and hence
\begin{equation}
\label{positano}
M(n_1,\dots,n_r)\ge 1{\rm ~if~}r\ge 2{\rm ~and~}n_1,n_2\ge 1.
\end{equation}
Case 1. $r=1$ and $n_1\ge 3$. We have $M(n_1,0)=0$ and $M(n_1,1)=1$, so 
we may assume that $m=1$. The word
$1^{n_1-1}2^{m}12$ 
is counted by $M(n_1,m+1)$, but not counted by the
$M(n_1,m)$ words amongst the $M(n_1,m+1)$ that come from  Lyndon extension
(since $1^{n_1-1}2^m1$ is not a Lyndon word). It follows that
$M(n_1,m+1)>M(n_1,m)$.\\
Case 2. $r\ge 2$, $n_1,n_2\ge 1$. Consider a Lyndon word counted
by $M(n_1,n_2-1,n_3,\dots,n_r,m+1)$ (such a word exists by (\ref{positano})) and
extend it with a 2. Since a Lyndon word counted by $M(n_1,n_2-1,n_3,\dots,n_r,m+1)$ starts
with a 1, this  will yield, by Proposition \ref{Lyndon}, again a Lyndon word
(counted by $M(n_1,\dots,n_r,m+1)$). This word is not amongst those coming
from Lyndon extension (they all end with a number $\ge 3$).\\
\indent On combining these results with Lemma \ref{allebegin} the proof of
part 1 is easily completed.\\
2) The argument of case 1 applies here as well. In addition we have to check now that
none of the words after extension is of the form $wwz$ (so as to avoid that they
are being counted as well under the $M(n_1/2,\dots,n_r/2,m/2)$ words that were
injected into $M(n_1,\dots,n_r,m+1)$ in the proof of Proposition \ref{mooizeg}.
A Lyndon word that is being extended is clearly not of the form $ww$. The only 
non-Lyndon word used in the previous argument, $1^{n_1-1}2^{m}12$ (with 
$m\ge 2$), is also not
of this form. \qed

\section{The proof of Theorem \ref{main1}}
Recall that if $f$ and $g$ are arithmetic functions, 
the classical M\"obius inversion
formula states that $g(n)=\sum_{d|n}f(d)$ iff
$f(n)=\sum_{d|n}\mu(d)g(n/d)$. Lemma \ref{inversie} is an analogue of
this result for sequences of formal series.
By writing $A^{(r)}(x)=\sum a_{j,r}x^j$ and $B^{(r)}(x)=\sum b_{j,r}x^j$, a
known M\"obius inversion formula for arithmetic functions in two variables is
obtained (see e.g. \cite{H}). Lemma \ref{inversie} is the main ingredient
in the proof of Theorem \ref{main1}.
\begin{Lem}
\label{inversie}
Let $\{A^{(r)}(z)\}_{r=1}^{\infty}$ and $\{B^{(r)}(z)\}_{r=1}^{\infty}$ be
two sequences of
formal series. Then
$$A^{(r)}(z)=\sum_{d|r}B^{(r/d)}(z^d){\rm ~iff~}
B^{(r)}(z)=\sum_{d|r}\mu(d)A^{(r/d)}(z^d).$$
\end{Lem}
{\it Proof}. We have
\begin{eqnarray}
\sum_{d|r}\mu(d)A^{(r/d)}(z^d)&=&\sum_{d|r}\mu(d)\sum_{e|{r\over
d}}B^{(r/de)}(z^{de})\nonumber\cr
&=&\sum_{m|r}(\sum_{d|m}\mu(d))B^{(r/m)}(z^m)\nonumber\cr
&=&B^{(r)}(z).
\end{eqnarray}
Conversely,
\begin{eqnarray}
\sum_{d|r}B^{(r/d)}(z^d)&=&\sum_{d|r}\sum_{e|{r\over
d}}\mu(e)A^{(r/de)}(z^{de})\nonumber\cr
&=&\sum_{m|r}(\sum_{e|m}\mu(e))A^{(r/m)}(z^m)\nonumber\cr
&=&A^{(r)}(z).
\end{eqnarray}
In both strings of identities we used that $\sum_{d|m}\mu(d)=0$ if $m>1$.
\qed\\

\noindent {\it Proof of Theorem } \ref{main1}.\\
1) Trivial.\\
2) Immediate from Lemma \ref{inversie} and the definition of ${\cal W}_f^{(r)}$.\\
3) Similar to the proof of Lemma \ref{minnetje}.\\
4) By Lemma \ref{inversie} it is enough to show that
$$\sum_{d|r}{r\over d}{\cal W}_{fg}^{({r\over d})}(z^d)=
\sum_{d|r}{r\over d}\sum_{[i,j]={r\over d}}(i,j)
{\cal W}_f^{(i)}(z^{r\over i}){\cal W}_g^{(j)}(z^{r\over j}).$$
By part 2 the left hand side of this purported identity equals $f(z)^rg(z)^r$.
The summation conditions on the right are equivalent to $i|r$ and $j|r$, that is,
$i$ and $j$ independently range over the divisors of $r$. Thus, noting that
$[i,j](i,j)=ij$, we obtain, for the right hand side,
$$\sum_{i|r}i{\cal W}_f^{(i)}(z^{r\over i})\sum_{j|r}j{\cal W}_g^{(j)}(z^{r\over j})
=f(z)^rg(z)^r,$$
where part 2 was used again.\\
5) By Lemma \ref{inversie} it is enough to show that
$$\sum_{d|r}{r\over d}{\cal W}_{f^k}^{({r\over d})}(z^d)=
\sum_{d|r}{r\over d}\sum_{[j,k]={rk\over d}}{jd\over r}
{\cal W}_f^{(j)}(z^{rk\over j}).$$
The left hand side is seen to equal $f(z)^{kr}$, the right hand side simplifies to
$\sum_{j|rk}j{\cal W}_f^{(j)}(z^{rk\over j})=f(z)^{rk}$.\\
6) Combine the identities of part 5 and part 6 (cf. the proof of Theorems 5
and 6 of \cite{MR}). \qed

\section{The proof of Theorems \ref{main2} and \ref{main3}}
We now have the necessary ingredients to prove Theorem \ref{main2}.\\
{\it Proof of Theorem} \ref{main2}.\\
1) Trivial.\\
2) Very similar to that of part 3: instead of (\ref{handig}) use
(\ref{handig2}). An
alternative proof is discussed in the next section.\\
3) Write $f(z)=\sum_j a_jz^j$. Let $k\ge 0$
be an arbitrary integer. In order to prove that the coefficient of $z^k$
in ${\cal W}_f^{(r)}(z)$ is
non-negative it is enough to prove this with $f(z)$ replaced by
$f_k(z)=\sum_{j=0}^k a_jz^j$.
Let $r=f_k(1)$. We apply equation (\ref{handig}) with $z_1,\dots,z_{a_1}=1$,
$z_{a_1+1},\dots,z_{a_1+a_2}=z$ etc.. Thus for example we write $2+z+z^2$
as $z_1+z_2+z_3+z_4$ with
$z_1=1$, $z_2=1$, $z_3=z$ and $z_4=z^2$. Since
every number of the form $M(m_1,\dots,m_r)$ is an integer, the result then
follows from (\ref{handig}).\\
4) Follows on combining part 3 above with part 3 of Theorem \ref{main1}.\\
5) From (\ref{handig}) and the definition of $M(n_1,\dots,n_r)$ (given
in (\ref{basaal})), we deduce that
\begin{eqnarray}
{1\over n}\sum_{d|n}\mu(d)(z_1^d+\dots+z_r^d)^{n\over d}&=&
{1\over n}\sum_{d|n}\mu(d)(z_1^d+\dots+z_{r-1}^d)^{n\over d}\nonumber\cr
&+&\sum_{n_1+\dots+n_r=n\atop n_j\ge 0,~n_r\ge 1}
M(n_1,\dots,n_r)z_1^{n_1}\cdots z_r^{n_r}.\cr
\end{eqnarray}
The remainder of the argument should be obvious to the reader (cf. the proof
of part 3). \qed\\

\noindent {\tt Remark 4}. Part 2 can also be proved using the theory of formal groups \cite{S}.\\

\noindent We now come to the proof of Theorem \ref{main3}.\\
{\it Proof of Theorem} \ref{main3}. 1) Write $f(z)=\sum_j a_jz^j$.
Put $w_1=\dots=w_{a_0}=0$,
$w_{a_0+1}=\dots=w_{a_0+a_1}=1$,
$w_{a_0+a_1+1}=\dots=w_{a_0+a_1+a_2}=2$, et cetera. On applying
(\ref{handig}), cf. the proof
of part 5, we see that the coefficient of $z^k$ in ${\cal W}_f^{(r)}(z)$
equals
\begin{equation}
\label{coefficient}
m_f(k,r)=\sum_{\sum_j n_jw_j=k,~\sum_j n_j=r\atop n_j\ge
0}M(\{n_j\}_{j=1}^{\infty}).
\end{equation}
The assumption $a_0>0$ implies that $w_1=0$. With each solution
of the system
$\sum_j n_jw_j=k$ and $\sum_j n_j=r$ we associate
$\{n_j'\}_{j=1}^{\infty}$ with
$n'_j=n_j$ for all $j\ge 2$ other and  with $n'_1=n_1+1$.
Note that $\sum_j n'_jw_j=k$ and $\sum_j n'_j=r+1$.
Note that the assignment $\{n_j\}_{j=1}^{\infty}\rightarrow \{n'_j\}_{j=1}^{\infty}$
is injective.
Using part 1 of Proposition \ref{mooizeg} we infer that
$M(\{n_j\}_{j=1}^{\infty})\le M(\{n'_j\}_{j=1}^{\infty})$.
This, in combination with (\ref{coefficient}) and the injectivity, yields that
$$m_f(k,r)\le
\sum_{\sum_j n'_jw_j=k\atop {\sum_j n'_j=r+1\atop n'_j\ge
0}}M(\{n'_j\}_{j=1}^{\infty})\le
\sum_{\sum_j m_jw_j=k\atop {\sum_j m_j=r+1\atop m_j\ge
0}}M(\{m_j\}_{j=1}^{\infty})=m_f(k,r+1).$$
2) The analogue of (\ref{handig}) reads
$$
{(-1)^n\over n}\sum_{d|n}\mu(d)(-z_1^d-\cdots-z_r^d)^{n\over
d}=\sum_{n_1+\cdots+n_r=n\atop n_j\ge 0}V_r(n_1,\cdots,n_r)z_1^{n_1}
\cdots z_r^{n_r},
$$
whence, cf. the proof of part 1, 
$$(-1)^rm_{-f}(k,r)=\sum_{\sum_j n_jw_j=k,~\sum_j n_j=r\atop n_j\ge
0}V_{\infty}(\{n_j\}_{j=1}^{\infty}),$$
where, if $t$ is such that $n_k=0$ for every $k\ge t+1$,
$V_{\infty}(\{n_j\}_{j=1}^{\infty})$ is defined as $V_t(n_1,\dots,n_t)$.
The remainder of the proof is quite similar to that of part 1.\\
3) Let $r\ge 1$. Suppose that $\sum_j n_jw_j=k$ and $\sum_j n_j=r$ with $n_j\ge 0$. Let
$t$ be the (unique) integer
such that $n_t\ge 1$ and $n_{t+j}=0$ for every 
$j\ge 1$. 
Let $t_1=t+a_{w_t}$. Since $a_{w_t}\le a_{w_t+1}$ (by assumption), it follows that
$w_{t_1}=w_t+1$.
To this solution
$\{n_j\}_{j=1}^{\infty}$ we associate $\{n'_j\}_{j=1}^{\infty}$ with
$n'_t=n_t-1$, $n'_{t_1}=1$
and $n'_j=n_j$ for $j\ne t,t_1$.  
Note that $\sum_j n'_jw_j=k+1$ and $\sum_j n'_j=r$.\\
3a) Let $t_2$ be such that $w_{t_2}=k$ 
(by assumption $f(z)\in \mathbb Z_{\ge 1}[[z]]$ and
hence $\{w_1,w_2,\dots\}=\mathbb Z_{\ge 0}$ and such a 
$t_2$ exists). Since $k\ge 1$ (by assumption), $t_2\ge 2$. 
Let $n_1=r-1$ and $n_{t_2}=1$ and set $n_j=0$ for $j\ne 1,t_2$.
Then $M(\{n_j\}_{j=1}^{\infty})\ge 1$ (by (\ref{positano}) and $M(0,1)=1$).
The result now follows by (\ref{coefficient}).\\
3b) From part
1 of Proposition \ref{mooizeg} and the fact that $M(n_1,\dots,n_r)$ is totally
symmetric in $n_1,\dots,n_r$, we infer that 
$M(\{n_j\}_{j=1}^{\infty})\le M(\{n'_j\}_{j=1}^{\infty})$.  
Using this, (\ref{coefficient}), and the injectivity of the
assignment $\{n_j\}\rightarrow \{n`_j\}$, we then obtain
\begin{equation}
\label{omhoogsammie}
m_f(k,r)\le
\sum_{\sum_j n'_jw_j=k+1\atop {\sum_j n'_j=r\atop n'_j\ge
0}}M(\{n'_j\}_{j=1}^{\infty})\le
\sum_{\sum_j m_jw_j=k+1\atop {\sum_j m_j=r\atop m_j\ge
0}}M(\{m_j\}_{j=1}^{\infty})=m_f(k+1,r).
\end{equation}
3c) If $r\ge 3$ and $k\ge 2$ one checks that the system of equations
$\sum_j n_jw_j=k$ and $\sum_j n_j=r$ has a solution in non-negative integers $n_j$
such that $n_t\ge 2$. For such a solution we have, by part 2 of Proposition
\ref{mooizeg}, $M(\{n_j\}_{j=1}^{\infty})<M(\{n'_j\}_{j=1}^{\infty})$. 
This ensures that the first inequality in (\ref{omhoogsammie}) becomes strict.\\
4) A variation of the proof of part 3 that is left to the reader. \qed\\

\noindent {\it The proof of Proposition} \ref{wassendwater}. Using (\ref{coefficient})
one finds that, for $c\ge 1$,
$$M(c;r)=\sum_{n_1+\dots +n_c=r\atop n_j\ge 0}M(n_1,\dots,n_c).$$
On invoking Proposition \ref{mooizeg} and Theorem \ref{atlonglast} the
result then follows after some calculation. \qed

\section{An alternative proof of part 2 of Theorem \ref{main2}}
An alternative proof of part 2 of Theorem \ref{main2} is obtained on combining
Theorem \ref{een} with Lemma \ref{lemmaeen}. Both the proof 
of Theorem \ref{een} and Lemma \ref{lemmaeen} are taken
from \cite{PNew} and
repeated here for the convenience of the reader.\\

\noindent {\it Proof of Theorem} \ref{een}. Using part 2 of Theorem
\ref{main1}
we deduce that
$$f(z)^r=\sum_{d|r}{r\over d}{\cal W}^{({r\over
d})}_f(z^d)=\sum_{d|r}{r\over d}
\sum_{j=0}^{\infty}m_f(j,{r\over d})z^{jd},$$
from which it is inferred that
$$\sum_{r=1}^{\infty}y^rf(z)^r=\sum_{j=0}^{\infty}\sum_{k=1}^{\infty}m_f(j,k)k\sum_{d=1}^{\infty}z^{jd}y^{kd}.$$
The latter identity with both sides divided out by $y$ can be rewritten as
$${f(z)\over
1-yf(z)}=\sum_{j=0}^{\infty}\sum_{k=1}^{\infty}{m_f(j,k)kz^jy^{k-1}\over
1-z^jy^k}.$$
Formal integration of both sides with respect to $y$ gives
$$-\log(1-yf(z))=-\sum_{j=0}^{\infty}\sum_{k=1}^{\infty}m_f(j,k)\log(1-z^jy^k),$$
whence the result follows. \qed\\

\noindent An {\it unital series} is, by definition, a series 
whose constant term equals 1. 
\begin{Prop} {\rm (Metropolis and Rota \cite[Proposition 1]{MR})}.
\label{uniec}
If $f\in \mathbb Z[[x]]$ is unital, then it has an unique expansion of
the form $f(z)=\prod_{n=1}^{\infty}(1-z^n)^{-e_n}$, where the $e_n$ 
are integers 
\end{Prop}
The following result generalizes this to two variables.
\begin{Lem} {\rm (Moree \cite{PNew})}.
\label{lemmaeen}
Suppose that $f(z,y)=\sum_{j,k}\alpha(j,k)z^jy^k$ where the $\alpha(j,k)$
are integers and $f(0,0)=0$. Then
there are unique integers $e(j,k)$ such that, as formal series, we have
$$1+f(z,y)=\prod_{j=0}^{\infty}\prod_{k=0\atop (j,k)\ne
(0,0)}^{\infty}(1-z^jy^k)^{e(j,k)}.$$
\end{Lem}
{\it Proof}. We say that $z^{j_1}y^{k_1}$ is of lower weight than
$z^{j_2}y^{k_2}$ if
$k_1<k_2$ or $k_1=k_2$ and
$j_1<j_2$. Suppose that $z^jy^k$ is the term of lowest weight appearing in
$f(z,y)$. Then consider
$(1+f(z,y))(1-z^jy^k)^{\alpha(j,k)}$. This can be rewritten as $1+g(z,y)$
where all the
coefficients of $g(z,y)$ are
integers and the term of lowest weight in $g(z,y)$ has strictly larger
weight than the term of lowest weight in
$f(z,y)$. Now iterate.\\
\indent It is not obvious from this argument that if we start with a
different weight ordering of the terms $X^jY^k$ we
end up with the same integers $e(j,k)$. Suppose that $h(X)$ has integer
coefficients, then the coefficients $e(n)$
in $1+h(X)=\prod_{n=1}^{\infty}(1-X^n)^{e(n)}$ are unique 
by Proposition \ref{uniec}. Hence, by setting $X=0$, respectively $Y=0$, we obtain
that $e(0,k)$, respectively $e(j,0)$ are uniquely determined. Setting
$Y=X^m$ we obtain
that $1+f(X,X^m)=\prod_{n=1}^{\infty}(1-X^n)^{f(n)}$, where $f(n)$ is
uniquely determined and
$f(2m)=e(2m,0)+e(m,1)+e(0,2)$. The uniqueness of $e(0,2)$, $e(2m,0)$ and
$f(2m)$ then implies the uniqueness of
$e(m,1)$. We
will complete the
proof by using induction. So suppose we have established that $e(j,k)$
with $k\le r$ for some $r\ge 1$ are uniquely determined.
Using that $f((r+2)m)=\sum_{k=0}^{r+2}e((r+2-k)m,k)$, we infer by
the induction hypothesis and using that $e(0,r+2)$ and $f((r+2)m)$ are
uniquely
determined, that $e(m,r+1)$ is uniquely determined. \qed

\section{Examples}
{\tt Example 1}. It turns out that the coefficients of the M\"obius transform in case
$f(z)=1+z$ have
several interesting properties. Note that, with $f(z)=1+z$, we have
\begin{equation}
\label{specialcase}
{\cal W}^{(r)}_f(z)={1\over r}\sum_{d|r}\mu(d)(1+z^d)^{r\over
d}=\sum_{j=0}^r m_f(j,r)z^j,
\end{equation}
where
$$m_f(j,r)={1\over r}\sum_{d|{\rm gcd}(j,r)}\mu(d)\left({{r\over d}\atop
{j\over d},{r-j\over d}}\right)
=M(j,r-j)\in \mathbb Z_{\ge 0}.$$
The numbers $M(j,r-j)$ also arise in the theory of relative partitions, see 
Remark 1. Witt's work \cite{W} yields the following result (where the formulation from
Proposition 2.10 of \cite{KK} with
$r=1$ is being used).
\begin{Prop}
\label{Propdrie}
Let $V=\bigoplus _{i,j=1}^{\infty } V_{(i,j)}$
be a $(\mathbb{Z}_{>0} \times \mathbb{Z}_{>0})$-graded vector space
over ${\mathbb C}$ with  $\mathrm{dim} V_{(i,j)}=1\in \mathbb{Z}_{>0}$
for all $i,j\ge 1$, and
let $L=\bigoplus _{m,n=1}^{\infty } L_{(m,n)}$ be the free Lie algebra
generated by $V$. We have
$${1\over n}\sum_{d|n}\mu(d)(1+z^d)^{n/d}=\sum_{j=0}^n {\rm
dim}(L_{(j,n-j)})z^j=\sum_{j=0}^n m_f(j,n)z^j.$$
\end{Prop}
Since $1+z$ is self-reciprocal, so is ${\cal W}_f^{(r)}(z)$ by part 1 of
Theorem \ref{main2}. It follows that dim$(L_{(j,n-j)})={\rm dim}(L_{(n-j,j)})$.\\

\noindent {\tt Example 2}.
Many constants in number theory have the form 
$\prod_{p>p_m}h(1/p)$, where $h$ is 
a rational function  and 
$h(z)=1+O(z^2)$ (as $z$ tends to zero) and the product is over all primes $p>p_m$, with $p_m$ 
the $m$th prime. Examples 
are the twin prime constant $T$ and $A$ the Artin
constant (defined in 
(\ref{arieartin})). We have the
formal identity $h(z)=\prod_{n=2}^{\infty}(1-z^n)^{-e_n}$, with the
$e_n$ uniquely determined integers (by 
Proposition \ref{uniec}). This identity can be used to expand $\prod_{p>p_m}h(1/p)$
in terms of the partial zeta function $\zeta_m(s)=\prod_{p\le p_m}(1-p^{-s})\zeta(s)$, 
where $\zeta(s)$ denotes the Riemann zeta function.
Formally we have
\begin{equation}
\label{constant}
\prod_{p>p_m}h({1\over p})=\prod_{p>p_m}
\prod_{n=2}^{\infty}(1-p^{-n})^{-e_n}=\prod_{n=2}^{\infty}\prod_{p>p_m}(1-p^{-n})^{-e_n}
=\prod_{n=2}^{\infty}\zeta_m(n)^{e_n}.
\end{equation}
For $m$ large enough it can be shown that such an identity always holds, see
Theorem 1 of \cite{Moreeconstant}. These identities can be used to evaluate constants
of this format with high numerical accuracy.\\
\indent To conclude we give a result concerning a class of more complicated constants
in which the Witt transform arises. These are the constants arising in the
left hand side of (\ref{ggg}), where $\chi$ is any Dirichlet character.  
\begin{Thm} {\rm (Moree \cite{PNew})}.
\label{algemener}
Suppose that $f(z)=\sum_{j\ge 1}a(j)z^j\in \mathbb Z[[z]]$. 
Let $j_0\ge 1$ denote the smallest 
integer such that $a(j_0)\ne 0$. Let $g(z)=\sum_{j\ge 1}|a(j)|z^j$.
Let $m_f(j,r)$ be defined as in Definition {\rm 1}.
Then, as formal power series in $y$ and $z$, one has
\begin{equation}
\label{dubbelproduct}
1-yf(z)=\prod_{k=1}^{\infty}\prod_{j=kj_0}^{\infty}(1-z^jy^k)^{m_f(j,k)},
\end{equation}
Moreover, the numbers $m_f(j,k)$ are integers.\\
\indent Let $\epsilon>0$ be 
fixed. The identity {\rm (\ref{dubbelproduct})} holds for all complex numbers $y$ and $z$
with $g(|z|)y<1-\epsilon$ and $|z|<\rho_c$, where $\rho_c$ is the radius of convergence of the Taylor
series of $g$ around $z=0$. If, moreover, $\rho_c>1/2$, $g(1/2)<1$ and
$\sum_p g({1\over p})$ converges, then
\begin{equation}
\label{ggg}
\prod_p\left(1-\chi(p)f({1\over p})\right)
= \prod_{k=1}^{\infty}\prod_{j=kj_0}^{\infty}L(j,\chi^k)^{-m_f(j,k)}.
\end{equation}
In the latter sum and product $p$ runs over all primes.
\end{Thm}
Recall that the Dirichlet L-series for $\chi^k$, $L(s,\chi^k)$, is defined, for Re$(s)>1$,
by $\sum_{n=1}^{\infty}\chi^k(n)n^{-s}$. Since Dirichlet L-series in integer values are
very easily evaluated with high decimal precision, this result allows one to evaluate 
with high decimal precision the constant
appearing on the left hand side of (\ref{ggg}). 
In the case of
the constants 
$$B_{\chi}=\prod_p\left(1+{[\chi(p)-1]p\over [p^2-\chi(p)](p-1)}\right),$$
arising in the study of some problems involving the multiplicative order, e.g. 
\cite{PNew, M, Pa},
one obtains from Theorem \ref{algemener} the following proposition:
\begin{Prop} {\rm (Moree {\rm \cite{PNew}})}.
\label{laatste}
Let $f(z)=-(1-z-z^2)^{-1}$.
We have
$$B_{\chi}=A{L(2,\chi)L(3,\chi)\over L(6,\chi^2)}
\prod_{r=1}^{\infty}\prod_{j=3r+1}^{\infty}L(j,\chi^r)^{-m_f(j-3r,r)},$$
where 
\begin{equation}
\label{arieartin}
A=\prod_p \left(1-{1\over p(p-1)}\right)=0.3739558136\cdots
\end{equation}
denotes the Artin constant and $(-1)^{r-1}e(j,r)=-m_f(j-3r,r)$.
\end{Prop} 
As a formal
series we have $1/(1-z-z^2)=\sum_{j\ge 0}F_{j+1}z^j$, with $F_j$ the $j$th Fibonacci number.
The Taylor coefficients of $(1-z-z^2)^{-r}$ are known as 
{\it convolved Fibonacci numbers}. Thus the numbers $e(j,r)$ are closely related to convolved Fibonacci numbers.
Numerical computation suggests that actually $e(j,r)\ge 1$ and, moreover, that these numbers enjoy certain
monotonocity properties in both the $j$ and $r$ direction. On using 
Theorem \ref{main2} various
of these numerical observations can be actually proved, see \cite{PFibo} for
details.\\

\noindent {\bf Acknowledgement}. There are too many mathematicians to thank here
for their kind e-mail comments on my questions. I'd like to especially
mention Profs Bryant, Fulman, Haukkanen, Jibladze, Petrogradsky, Reiner, Ruskey, Stanton 
and Stienstra.
Prof. Ruskey's 
suggestion of using Lyndon words turned out to be very useful. Profs. Feingold, 
Kang and
Mays  kindly sent me
some reprints. I thank Dr. Christa Binder 
for sending me some biographical material
concerning the Austrian mathematician R. Daublebsky von Sterneck (1871-1928) and
several of his illustrous relatives.

\medskip\noindent {\footnotesize Max-Planck-Institut,
Vivatsgasse 7, D-53111 Bonn, Germany.\\
e-mail: moree@mpim-bonn.mpg.de }

\end{document}